\documentclass{article}
 
\usepackage{amsmath}
\usepackage{amssymb}
\usepackage{amscd}
\usepackage[all]{xy}

          \newtheorem{df}{Definition}[section]
          \newtheorem{pr}[df]{Proposition}
          \newtheorem{theorem}[df]{Theorem}
          \newtheorem{lem}[df]{Lemma}
          \newtheorem{cor}[df]{Corollary}
          \newtheorem{rem}[df]{\it Remark}

          \newtheorem{prob}[df]{Problem}
          \newtheorem{exam}[df]{Example}

          \newcommand{\qed}{~$\Box$\newline\hskip 0.18cm}
          \newcommand{\proof}{
                             \noindent{\bf Proof.~~}}
          \newcommand{\dslash}{\slash\hspace{-0.1cm}\slash}

            \makeatletter
            \@addtoreset{equation}{section}
            \makeatother
         
         \newcommand{\mapright}[1]{%
           \smash{\mathop{%
           \hbox to 1cm{\rightarrowfill}}\limits^{#1}}}
        
         \newcommand{\mapleft}[1]{%
           \smash{\mathop{%
           \hbox to 1cm{\leftarrowfill}}\limits_{#1}}}
           
         \newcommand{\shortmapright}[1]{%
           \smash{\mathop{%
           \hbox to .7cm{\rightarrowfill}}\limits^{#1}}}
        
         \newcommand{\shortmapleft}[1]{%
           \smash{\mathop{%
           \hbox to .7cm{\leftarrowfill}}\limits_{#1}}}
           
\begin{document}

       
       \title{Valuative characterization of central extensions \\ of algebraic tori on Krull domains
       \footnote{This is a version on  August 3, 2016.}}
       
         
         \author{{Haruhisa} {\sc Nakajima \footnote{Partially supported by  
        Grant No. 26400019:  the Japan Society for the Promotion of Sciences.}}\\
        \small { Department of Mathematics}, 
        \small {\sc  J. F. Oberlin University}\\
        \small Tokiwa-machi, {\sc Machida}, {\sc Tokyo}  194-0294,
         {\sc JAPAN}} 
         \date{ }
         \maketitle

        
         \begin{abstract}
            
            Let $G$ be an affine algebraic group with 
 an algebraic torus $G^{0}$ over an algebraically closed field $K$ of
an arbitrary  characteristic $p$.   We show  a
criterion  for $G$  to be a finite central extension of $G^{0}$
in terms of invariant theory of  all  regular actions of any closed subgroup $H$ containing 
$Z_{G}(G^{0})$ on affine
 Krull $K$-schemes such that  invariant rational functions are locally
fractions of invariant regular functions.
 Consider an affine Krull  $H$-scheme $X =
{\rm Spec} (R)$ and a prime ideal ${\mathfrak P}$ of $R$ with ${\rm ht} ({\mathfrak P})
= {\rm ht} ({\mathfrak P}\cap R^{H}) = 1$. 
 Let ${\mathcal I}({\mathfrak P})$
denote the inertia group of ${\mathfrak P}$ under the action of $H$. 
The group  $G$ is 
central over $G^{0}$ if and only if the fraction
${\rm e} ({\mathfrak P} , {\mathfrak P}\cap R^{H}) / {\rm e} ( {\mathfrak P}, 
{\mathfrak P}\cap R^{{\mathcal I}({\mathfrak P})})$
of ramification indices is equal to $1$  ($p=0$)  or to  the $p$-part of 
the order  of the group of weights of $G^{0}$ on $R^{{\mathcal I}({\mathfrak P}))}$
vanishing on $R^{{\mathcal I}({\mathfrak P}))}/{\mathfrak P}\cap R^{{\mathcal I}({\mathfrak P})}$ ($p>0$) for 
an arbitrary $X$ and  ${\mathfrak P}$.

              
            \end{abstract}

    \normalsize
         \small \section{Introduction}\normalsize
 \paragraph{\bf 1.A.}   We consider  affine algebraic groups and affine schemes over  a fixed algebraically closed field $K$ of  an arbitrary characteristic $p$.  For an affine group $G$,
 denote by $G^{0}$ its identity component.  
 Let $\rho:  G \to GL(V(\rho))$ be a finite dimensional rational representation of $G$
 and $A = K[V({\rho})]$ a $K$-algebra of polynomial functions on $V(\rho)$.
 We have the inclusions 
\begin{equation} A^{G}  \subseteq  A^{G^{0}} \subseteq A^{(G^{0})'}\subseteq A^{{\mathfrak Rad}_{u}(G^{0})}\end{equation}
of rings of invariants,  where ${\mathfrak Rad}_{u}(G^{0})$ is the unipotent radical of $G^{0}$
 and $(G^{0})'$ is the inverse image in $G^{0}$ of the semi-simple part of 
 $G^{0}/{\mathfrak Rad}_{u}(G^{0})$.  These rings are Krull domains but, 
 they are not necessarily 
 finitely generated over $K$ if  ${\mathfrak Rad}_{u}(G^{0})$ is non-trivial. Since $A^{G}$ is obtained as  a ring of invariants of the  Krull $K$-domain 
 $A^{(G^{0})'}$ (resp. $A^{{\mathfrak R}(G^{0})}$) 
 under the  action of ${G/(G^{0})'}$ (resp. $G/{{\mathfrak Rad}_{u}(G^{0})}$) whose identity component  is an algebraic torus (resp. a reductive group). 
 On the other hand by slice \'etale theorem (e.g., \cite{Po}),  it needs to study on  invariants of 
 the stabilizer $G_{P}$ of a suitable point $P$, which may not be connected.  Moreover
 reductivity of algebraic groups is 
 characterized by pseudo-reflections of actions on Krull $K$-domains
  (cf. \cite{Nak5}). 
  Consequently
in order to examine  $A^{G}$ we have to study rings of invariants
of Krull $K$-domains with the above group actions in general settings.

 On  (1.1), in the case where  $G^{0}$ is semisimple and $p =0$,  
 if $A^{G}$ is polynomial ring then $A^{G^{0}}$
 is a complete intersection (cf.  \cite{Panyushev}). This  seems to be useful to classify
 such representations of $G$. The author
  studied  in \cite{Nak1} the representation 
   $\rho$ of $G$ with  an algebraic torus $G^{0}$
 such that $R^{G}$ is a polynomial
 ring in the case of $p= 0$.  Our result suggests {\it  there are differences
 in invariant theory  of 
finite central extensions  and  finite non-central extensions of a torus
related to the role of pseudo-reflections.}
Hence we will  clarify  the background of this phenomenon in this paper
and  advance  Problem 1.1 as follows. 
  
  \paragraph{\bf 1.B.}   
   Affine $K$-schemes $X$ are 
        affine schemes of commutative $K$-algebras  ${\mathcal O}(X)= R$ of
        global sections   which are not necessarily   finite generated as algebras over $K$. Let ${\mathcal Q}(R)$ denote the total quotient ring 
        of $R$. 
      We say an action 
        $(X, G)$ of  an affine algebraic group $G$ on an affine $K$-scheme $X$
         is {\it regular}, when $G$ acts {\it rationally}
        on the $K$-algebra $R$ as $K$-algebra automorphisms, i.e.,  $R$ is regarded as a
        rational $G$-module. In this case
        we denote by $X\dslash G$  the affine $K$-scheme defined by 
        $R^{G}$ which is the $K$-subalgebra consisting of all invariants of
         $R$ under the action of $G$ and by $$\pi _{X, G} : X \to X\dslash G$$
         the morphism  induced from the inclusion 
         $R^{G} \hookrightarrow R$. 
        Furthermore
        $(X, G)$ is said to be {\it effective} if $\rm{Ker} (G \to \rm{Aut} ~X)$ is finite. 
      
       \paragraph{\bf 1.C.}   
       The affine $K$-scheme $X = {\rm Spec}  (R)$ is said to be Krull, if $R$ is a Krull 
       $K$-domain. 
        In this case 
    let ${\rm e} ({\mathfrak P}, {\mathfrak P}\cap R^{G})$ be the ramification index is defined to satisfy
     $({\mathfrak P} R_{\mathfrak P})^{{\rm e} ({\mathfrak P}, {\mathfrak P}\cap R^{G}) }
     = (\mathfrak P\cap R^{G} )R_{\mathfrak P}$
     for a prime ideal ${\mathfrak P}$ of $R$ of ${\rm ht} ({\mathfrak P}) =
     {\rm ht} ({\mathfrak P}\cap R^{G}) =1$.  The ramification indices are related 
     to the order of pseudo-reflections of $G$ which are  important in the study 
     on  relative invariants  (cf. \cite{St1, Nak0, Nak5-0}). 
     
          \begin{prob}{\rm (Reduction of ramification indices to
          finite subgroup quotients)}  Is there a finite subgroup $I_{\mathfrak P}$ of $G$
          such that ${\rm e} ({\mathfrak P}, {\mathfrak P}\cap R^{G})$
          is equal to  a product of 
           ${\rm e} ({\mathfrak P}, {\mathfrak P}\cap R^{I_{\mathfrak P}})$
          and  $p^{n}$ for some computable $n \geqq 0$?
 
 \end{prob}
 This is affirmative  in the restricted case  treated in \cite{Nak1} and
 plays a key role in the classification of representations of disconnected tori
 with regular rings of invariants.  Now put
     $${\mathcal I}_{G} ({\mathfrak P}) := \{\sigma \in G \mid \sigma(x) -x \in {\mathfrak P}
     ~(\forall x \in R) \}.$$
    This group is {\it finite on} $X$ for $G$ with a reductive $G^{0}$ (cf. \cite{Nak5}).
      In this paper we solve Problem 1.1 completely for
        $G$ with an algebraic torus $G^{0}$. It should be noted that this problem is not
        affirmative,  if  $G^{0}$ is not a torus (cf. Example 6.2). 
        
         \paragraph{\bf 1.D.}   
     The present  paper is a continuation of \cite{Nak2}   and completes its  studies.  We will establish
        a  characterization of a finite central extension of a torus in terms of 
        ramification indices in invariant theory. In order to study on them, 
         it is natural to treat  regular actions on  {\it Krull} $K$-domains. 
          Invariants of finite central extensions of algebraic tori with relations to  pseudo-reflections
        (or ramification theory) 
       are  useful  in studying the obstruction groups for cofreeness defined in 
       \cite{Nak5-0}. Those groups evaluate failure  of the Russian conjecture 
        for actions of algebraic tori on normal varieties (cf. Sect.1 and Examples 5.6, 5.7
        of \cite{Nak5-0}) in characteristic zero.

\paragraph{{\bf 1.F.}}

 Let  ${\mathfrak X}(G)$ be the group of rational characters of $G$ expressed as
an additive group with zero denoted to $0_{G}$. For a rational $G$-module  $M$, put
\begin{equation*}\left\{
\begin{aligned}  &{\mathfrak X}(G)^{M} : = \{ \chi \in  {\mathfrak X}(G) \mid  M_{\chi}\not=\{0\}\}\\
&{\mathfrak X}(G)_{M} :=  \{ \chi \in  {\mathfrak X}(G)^{M} \mid  
M_{-\chi}\not=\{0\}
\}\end{aligned}\right..
\end{equation*}
 Here $M_{\chi}$ stands for the space
  $\{ x\in M \mid \sigma (x) =\chi (\sigma) x\}$ of {\it relative invariants} 
  of $G$ in $M$ relative to $\chi$. Clearly ${\mathfrak X}(G)_{M}$ is a subgroup
  of ${\mathfrak X}(G)$,  which shall be used in Sect. 3.

\paragraph{1.G.} Let us explain  our main result  as follows. 
Suppose that $G^{0}$ is an algebraic torus and let  $H$ be a closed subgroup
of $G$ containing the centralizer $Z_{G}(G^{0})$ of $G^{0}$ in $G$. 
Consider    a regular action $(X, H)$
of $H$ on an affine Krull  $K$-scheme $X = {\rm Spec} (R)$. 
 For a prime ideal  ${\mathfrak P}$  of $R$
such that  ${\rm ht} ({\mathfrak P} )= {\rm ht} ({\mathfrak P}\cap R^{H}) =1$, we define 
$$\Delta_{H, {\mathfrak P}} := \left\{\begin{array}{lr}
\text{$p$-part of the index} ~\left[\langle {\mathfrak X} (G^{0})^{ R^{{\mathcal I}_{H}({\mathfrak P}))}}\rangle : \langle{\mathfrak X}(G^0)^{R^{{\mathcal I}_{H}({\mathfrak P})}/{\mathfrak P}\cap R^{{\mathcal I}_{H}({\mathfrak P})}}\rangle \right]  & (p > 0)\\
1  & (p =0)\end{array}\right. $$
The main theorem of this paper  is as follows:
\begin{theorem}\label{mainresult}
Suppose that $G^0$ is an algebraic torus.  Then the following conditions are equivalent:
\begin{itemize}
\item[(i)]  $G = Z_G(G^0)$
\item[(ii)] For an arbitrary   closed subgroup $H$ of   $G$ containing $Z_G(G^0)$, 
the following 
conditions hold  for any 
effective regular action   $(X, H)$
on an arbitrary affine Krull $K$-scheme $X = {\rm Spec} (R)$ such that ${\mathcal Q}(R^{G^{0}}) = {\mathcal Q}(R)^{G^{0}}$:   If ${\mathfrak P}$ is any prime ideal of $R$
with ${\rm ht} ({\mathfrak P} )= {\rm ht} ({\mathfrak P}\cap R^{H}) =1$, 
the following equality holds: 
$${\rm e}({\mathfrak P}, {\mathfrak P}\cap R^H)= 
{\rm e}({\mathfrak P}, {\mathfrak P}\cap R^{{\mathcal I}_H({\mathfrak P})}) 
\cdot   \Delta_{H, {\mathfrak P}}.
$$

\end{itemize}
\end{theorem}

 It should be noted that  the right hand side of  the equalities in {\it (ii)} are
  obtained  from  datum on finite group quotients, because  ${\mathcal I}_H({\mathfrak P})$  is finite on $X$.  The implication {\it (i)} $\Rightarrow$ {\it (ii)} can be regarded as
  a generalization of  a known result of reduced ramification indices of Krull domains under the 
  actions of  finite Galois groups.

\paragraph{1.H.} We prepare some results on invariant rational functions in Sect.2 which are helpful for replacing stabilities of algebraic group actions with conditions on 
quotient fields. Also we  prepare the  results  on  rational characters in Sect.3
which are useful  in the calculation of ramification indices.  In Sect.4
the structure of inertia groups of toric quotients is studied. In Sect.5
some auxiliary results in Sect.3 of \cite{Nak2} are generalized  
from the case of normal varieties to one of  Krull schemes. 
The results mentioned in Sect.1 $\sim$ Sect.4  are devoted to the
proof the theorem as above which is recognized as a generalization of the main result of \cite{Nak2}.  Its proof shall be
completed in the last section of this paper. 
 \smallskip 
\paragraph{Notations and Comments} The following notations are used without explanation:\\ 
\vskip -0.3cm
 \noindent $\bullet$ Let ~$p$ denote always the characteristic of  an fixed 
  algebraically closed field $K$. \\
\vskip -0.3cm  
\noindent 
$\bullet$~Let ${\rm U}(R)$ (resp. ${\mathcal Q}(R)$) denote the unit group 
(resp.  the total quotient ring) of a commutative ring $R$. Let ${\rm ht} ({\mathfrak I})$
denote the height of an ideal ${\mathfrak I}$ of $R$. \\
\vskip -0.3cm
\noindent
 $\bullet$~For a morphism $\gamma : H \to G$ of groups, let $Z_G(H)$ denote the
centralizer of a subset 
$\gamma(H)$ in  $G$ (i.e., $Z_{G}(H) =\{ \sigma\in G \mid 
\sigma \cdot \gamma(\tau) = \gamma (\tau) \cdot \sigma ~(\tau \in H) \}$). \\
\vskip -0.3cm
\noindent
$\bullet$~For a subset $C$ of a group $G$, let $\langle C\rangle$
denote the subgroup of $G$ generated by $C$.  For a subgroup $N$
of $G$, let $[G : N]$ denote the index of $N$ in $G$. \\
\vskip -0.3cm
\noindent 
$\bullet$~In the case where a group $G$ acts on a set $X$, consider
the associated morphism $G \to \mathrm{Sym} (X)$ from $G$ to the
(symmetric) group of all bijective transformations on $X$. Let $\sigma\vert_{X}$
denote the image of $\sigma\in G$ under the
canonical morphism $G \to \mathrm{Sym} (X)$. Moreover
let $G\vert_{X}$ be $\mathrm{Im} (G \to \mathrm{Sym} (X))$. We call ${\rm Ker} (G \to  {\rm Aut}~X)$ the {\it ineffective kernel} of the action $(X, G)$. If it is trivial (resp. finite), the action of $G$ on 
$X$ is said to be {\it  faithful} (resp. {\it effective}).  \\
\vskip -0.3cm
\noindent $\bullet$ ~The $p$-part $n_{[p]}$ and the   $p'$-part $n_{[p']}$ of $n \in \mbox{\boldmath $N$}$ are the natural numbers in such a way that $n = p^{e}\cdot n_{[p']}$,  $e\geqq 0$, $n_{[p']}$ is not divisible by $p$ and $n_{[p]} =p^{e}$
if p is a prime number. When $p =0$, we define that the $0$-part $n_{[0]} = 1$ and the $0'$-part $n_{[0']} = n$ of $n \in \mbox{\boldmath $N$}$. \\
\vskip -0.3cm
\noindent 
$\bullet$~The symbol  $\sharp (Y)$  denotes the cardinality of a set $Y$. The symbol  $\sharp_{p} (Y)$ (resp.   $\sharp_{p'} (Y)$) denotes the
$p$-part (resp.  $p'$-part ) of $\sharp (Y)$ of   a finite set $Y$ if  $p$ is prime or 0. \\ 
\vskip -0.3cm
\noindent 
$\bullet$~For subsets $A$, $B$ of $C$, let $A\backslash B$ denote the
difference set. \\

 \section{Preliminaries}

\paragraph{2.A.}   For a commutative ring $A$, let 
${\rm Ht}_{1}(A)$ be the set consisting of prime ideals of $A$ of height $1$. 
 For a prime ideal ${\mathfrak p}\in {\rm Ht}_{1}(B)$ of a subring $B$ of $A$,  let 
$${\rm Over}_{\mathfrak p} (A) = \{
{\mathfrak P} \in {\rm Ht}_{1}(A) \mid  {\mathfrak P}\cap B = {\mathfrak p}\}$$
and ${\rm Ht}_{1}(A, B) = \{ {\mathfrak P} \in {\rm Ht}_{1}(A) \mid 
{\rm ht} ({\mathfrak P}\cap B) = 1 \}$. 

\paragraph{2.B.} Let $H$ be a group acting on  $A$ as ring automorphisms. For 
a prime ideal $\mathfrak Q$ of $A$, let 
$${\mathcal D}_{H}({\mathfrak Q}) = \{ \sigma \in H\mid \sigma({\mathfrak Q} ) = {\mathfrak Q}\},$$
$${\mathcal I}_{H}({\mathfrak Q}) = {\mathrm{Ker}} \left({\mathcal D}_{H} ({\mathfrak Q})\underset{\text can.}\rightarrow {\rm{ Aut}} (A/{\mathfrak Q})\right)$$
which are respectively called the decomposition (reps. inertia) group of $H$ at ${\mathfrak Q}$. Let $A^{H}$ denote the subring of $A$ consisting of invariant elements   of $A$
under the action of $H$, whose elements are called {\it invariants} of $A$ under $H$.  If $A$ is an integrally closed domain, then so is $A^{H}$. For convenience, we note that the Galois theory of integrally closed domains
with finite group actions  can  be found in Chap. V of  \cite{LR} and 
 its \'etale version is  given in Exp. V of \cite{Gr} with locally noetherian conditions. 
 The  ramification theory for Dedekind domains with finite group actions is 
 treated in  Chap. V of \cite{zariski}.  Clearly
 ${\mathcal Q} (A) ^{H}  = {\mathcal Q} (A^{H})$ for a finite subgroup $H$ of ${\rm Aut} A$
 of an integral domain $A$.  However  the equality ${\mathcal Q} (A) ^{H}  = {\mathcal Q} (A^{H})$
 is not true unless $H$ is finite as follows. 

\begin{exam}\label{stability}  \rm Let $A = K[X_{1}, X_{2}, X_{3}, X_{4}]$ be a polynomial ring over $K$ with four variables  $X_{i}~ (1\leq i\leq 4)$. Let 
$G$ be an algebraic  torus $(\mathbf{G_{m}})^{2}$ of rank $2$. 
Suppose  any $\sigma_{s, t}= (s, t) \in \mathbf{G_{m}}\times \mathbf{G_{m}}$
acts $A$ as $K$-algebra 
automorphisms
by $$\sigma_{s, t} \left(\left[\begin{array}{c} X_{1} \\ X_{2} \\ X_{3} \\ X_{4} 
\end{array}\right]\right)=
\left[\begin{array}{ccccc} s^{-1} & & & \\ & s &  & \\ & & st &   \\   & & & t
\end{array}\right]  \left[\begin{array}{c} X_{1} \\ X_{2} \\ X_{3} \\ X_{4} 
\end{array}\right].$$ Then ${\mathcal Q}(A^{G}) = K(X_{1}X_{2})
\subsetneq {\mathcal Q} (A)^{ G} = K\left(
X_{1}X_{2}, 
\dfrac{X_{1}X_{3}}{X_{4}}\right)$.

\end{exam}

 \begin{lem}\label{qf}  Let $G$ be a group acting on an integrally closed domain $A$ as ring automorphisms.
 Let $N$ be a finite normal subgroup of $G$ and $H$ a normal subgroup of $G$ of a finite index (i.e., $[G : H] < \infty$)
 such that $N \subseteq H$.  Then ${\mathcal Q}(A^G) = {\mathcal Q}(A)^G$ if and only if
 ${\mathcal Q}(A^H) = {\mathcal Q}(A^N)^{H/N}$. 
 \end{lem}
\proof 
The {\it if part} of the assertion follows  easily from finiteness of $G/H$.  Suppose that ${\mathcal Q}(A^G) = {\mathcal Q}(A)^G$. 
Let $\widetilde{A^G}$ be the integral closure of $A^G$ in ${\mathcal Q}(A)^H$. 
As ${\mathcal Q}(A)^{H}$ is a Galois extension over ${\mathcal Q}(A)^{G}$
under the action of $G/H$,  we easily see ${\mathcal Q}(\widetilde{A^G})
= {\mathcal Q}(A)^H$. Since any element  of  $\widetilde{A^G}$ is integral over $A$ and $A$ is integrally closed, 
we must have $$\widetilde{A^G}\subseteq A\cap {\mathcal Q}(A)^H = A^H.$$ which shows
${\mathcal Q}(A^H) = {\mathcal Q}(A)^H$.  Then we see ${\mathcal Q}(A^H) = 
{\mathcal Q}(A^{N})^{H/N}$.\qed

 \paragraph{2.C.}        
        An affine scheme $Z= {\rm Spec} (R)$ is said to  be   Krull if $ R$
        is a Krull domain.  Let 
        $\mathrm{v}_{R, \mathfrak P}$ denote
       the discrete valuation of $Z$ defined 
       by the prime ideal ${\mathfrak P} \in {\rm Ht}_{1}(R)$. 
  It is well known that, for  any subfield $F$ of  
        ${\mathcal Q}(R)$, the subring $R\cap F$ is also 
        Krull.  If $R\cap F$ is not a field, 
        the set ${\rm Over}_{\mathfrak p} (R)$  is known to be  non-empty for any $\mathfrak p \in {\rm Ht}_{1}(R\cap F)$
        by the independence theorem of 
valuations (cf. \cite{LR, Fossum, BourbakiC, Magid}).
          {\it The ramification index of ${\mathfrak P}$ over} ${\mathfrak p}$
               is defined to be   $${\rm e} ({\mathfrak P}, {\mathfrak p}) : = 
        {\rm v}_{  R, {\mathfrak P} }( {\mathfrak p}  \it R_{\mathfrak P})$$
        where ${\mathfrak P}\in {\rm Over}_{\mathfrak p} (R)$. The 
        $p$-part (resp. the  $p'$-part) (cf. Sect. 1) of ${\rm e} ({\mathfrak P}, {\mathfrak p})$ 
        stands for 
        ${\rm e}_{p} ({\mathfrak P}, {\mathfrak p})$ (resp. 
        ${\rm e}_{p'} ({\mathfrak P}, {\mathfrak p}$)).

\begin{lem}\label{quotient}  Let $R$ be a Krull domain acted by a group $G$ as ring automorphisms. Then $R^{G}$ is a Krull domain.  In the case where $R^{G}$ is not a field, 
the set ${\rm Over}_{\mathfrak p}(R)$ is non-empty for  any ${\mathfrak p}\in {\rm Ht}_{1} (R^{G})$ and   
 the
following conditions are equivalent: 
\begin{itemize}
\item[(i)] ${\mathcal Q}(R^G) = {\mathcal Q}(R)^G$.
\item[(ii)] $R_{\mathfrak P}\cap {\mathcal Q}(R)^G = R^G_{{\mathfrak P}\cap R^G}$
for  ${\mathfrak P}\in {\rm Ht}_1(R, R^G)$. 


\end{itemize} 

\end{lem}
\proof    As $R^{G} = R\cap {\mathcal Q} (R)^{G}$, the former two assertions follow from the preceding paragraph  (cf. \cite{Magid}). 
Since $R_{\mathfrak P}$ with a prime ideal $\mathfrak P$ of $R$ of height 1 is a valuation ring of ${\mathcal Q}(R)$, for any nonzero element $x$ of ${\mathcal Q}(R)^G$, $x$ or $x^{-1}$ belongs to $R_{\mathfrak P}\cap {\mathcal Q}(R)^G$, which  shows the implication {\it (ii)} $\Rightarrow {\it (i)}$. 

 {\it (i)} $\Rightarrow$ {\it (ii)} : Suppose that ${\mathcal Q}(R^G) = {\mathcal Q}(R)^G$. For 
 any ${\mathfrak P}\in {\rm Ht}_1(R, R^G)$, $R_{\mathfrak P}\cap {\mathcal Q}(R)^G$ is a
 valuation ring of ${\mathcal Q}(R^G)$ contains the discrete valuation ring $R^G_{{\mathfrak P}\cap R^G}$ and  by Chap. 4 of \cite{Matsumura} (e.g. by
 Chap. 7 of \cite{BourbakiC}) of  the ring  $R_{\mathfrak P}\cap {\mathcal Q}(R)^G$ is a localization of
 $R^G_{{\mathfrak P}\cap R^G}$ at its prime ideal.  Unless  {\it (ii)} holds, we have 
 $R_{\mathfrak P}\cap {\mathcal Q}(R)^G = {\mathcal Q}(R^G).$ Then as 
 $${\mathfrak P}R_{\mathfrak P} \supseteq {\mathfrak P}R_{\mathfrak P}\cap {\mathcal Q}(R)^G
 \supseteq {\mathfrak P}\cap R^G\not=\{0\},$$ ${\mathfrak P}R_{\mathfrak P}$ contains a unit,  which is a contradiction.  \qed


\begin{pr}\label{etale}
Let $A$ be an integrally closed domain acted faithfully
 by  a finite group $H$ as automorphisms and   $L$ a subgroup   of $H$.  Let $\mathfrak Q$ be a prime ideal of $A$ and suppose that $A^{H}_{ {\mathfrak  Q}\cap A^{H}}$
is noetherian. 
Then 
the canonical monomorphism  $$A^{H}_{{\mathfrak Q}\cap A^{H}} \to  
A^{L}_{{\mathfrak Q}\cap A^{L}}$$
is \'etale if and only if $L \supset {\mathcal I}_{H}({\mathfrak Q})$. 

\end{pr}
\proof  Since  $A^{H}_{{\mathfrak Q}\cap A^{H}}\otimes_{A^{H}} A$ is an integral closure
of the integrally closed domain $A^{H}_{{\mathfrak Q}\cap A^{H}}$
in the finite separable extension ${\mathcal Q}(A)$ over ${\mathcal Q}(A)^{H}$, 
we see that $A^{H}_{{\mathfrak Q}\cap A^{H}}\otimes_{A^{H}} A$ is  finitely generated
as an  $A^{H}_{{\mathfrak Q}\cap A^{H}}$-module  (e.g., Chap. V of  \cite{zariski}). 
Exchanging the localization $A^{H}_{{\mathfrak Q}\cap A^{H}}\otimes _{A^{H}} A$ with $A$, we may suppose
that the rings $A$, $A^{L}$ and $A^{H}$ are noetherian
and $A$ is finitely generated as an $A^{H}$-module. Then the assertion follows
from Proposition 2.2 and Corollary 2.4 of Expos\'e V of  \cite{Gr}. \qed
\begin{cor}\label{etale2} 
Let $A$ be a Krull domain acted faithfully by  a finite group $H$ as automorphisms and   $L$ a subgroup   of $H$.  For a prime ideal ${\mathfrak Q}\in {\rm Ht}_{1}(A)$, 
the canonical monomorphism  $$A^{H}_{{\mathfrak Q}\cap A^{H}} \to  
A^{L}_{{\mathfrak Q}\cap A^{L}}$$
is \'etale if and only if $L \supset {\mathcal I}_{H}({\mathfrak Q})$. 
\end{cor}
\proof 
Since  $A^{H}_{{\mathfrak Q}\cap A^{H}}$ is a discrete valuation ring, this assertion  follows
from Proposition \ref{etale}.   \qed 

\paragraph{2.D.} The pseudo-reflections of finite linear groups and its slight modification are defined and studied in \cite{Bourbaki, St1, Nak0}.  We generalize  the
group generated by  pseudo-reflections
as follows: 
For  a Krull domain $R$ acted by a group $G$ as ring automorphisms,
let $${\mathfrak R}(R, G) = \left\langle  
\bigcup _{{\mathfrak P}\in {\rm Ht}_1(R, R^G)} {\mathcal I}_G({\mathfrak P})
\right\rangle \subseteq G$$ which is called
{\it the pseudo-reflection group of the action of $G$ on $R$}. 

In the case where {\it $G\vert _{R}$ is finite,}
for any ${\mathfrak q} \in {\mathrm Ht}_{1}(R^{G})$, consider  $R^{G}_{\mathfrak q}\otimes_{R^{G}} R$ with the action of $G$. Then
we can apply the ramification theory of Dedekind domains (e.g.,  Chap. V of \cite{zariski})
to this action and must have
${\rm e} ({\mathfrak Q}, {\mathfrak Q}\cap R^{G}) = 
{\rm e} ({\mathfrak Q}, {\mathfrak Q}\cap R^{{\mathcal I}_{G}({\mathfrak Q})})$ and 
\begin{equation}\label{eq1} \mathrm{e}_{p'} ({\mathfrak Q}, {\mathfrak Q}\cap R^{{\mathcal I}_{G}({\mathfrak Q})})
	= \sharp_{p'}({\mathcal I}_{G}({\mathfrak Q})\vert _{R}) \end{equation}
for ${\mathfrak Q} \in {\rm Over}_{{\mathfrak q}} (R)$, 	where $p$ denotes  the characteristic of $R$.  The last equality holds in the case where ${\mathcal I}_{G}({\mathfrak Q})\vert _{R}$
is finite, even if $G$ is not  finite on $R$, because  ${\mathcal I}_{H}({\mathfrak Q}) = H$
for $H = {\mathcal I}_{G}({\mathfrak Q})$. 

 \small \section{Character Groups and Inertial Quotients}\normalsize
 
          {\it  Hereafter in this paper we suppose that  $G$ is an affine algebraic group 
          over $K$}.  A regular action $(X, G)$ on an affine variety $X$
           is said to be {\it stable}, 
           if $X$ contains a non-empty
open subset consisting of closed $G^{0}$-orbits (or, equivalently, closed $G$-orbits)  (cf. \cite{Po})
           
\begin{rem}\label{stabilitycriteria} Suppose that  $G^{0}$ is an algebraic torus 
and $(X, G)$ is  a regular action of $G$ on an affine integral $K$-scheme $X ={\rm Spec} (R)$.
If \begin{equation}\label{range} {\mathfrak X}(G^{0})^{R}= 
{\mathfrak X}(G^{0})_{R}\end{equation}
(for the notation, see {\bf 1.F.} of Sect.1),  then ${\mathcal Q} (R^{G}) = {\mathcal Q} (R)^{G}$. 
 Especially in the case where  $R$ is finitely generated over $K$ as a $K$-algebra, $(X, G)$ (or $(X, G^{0})$) is stable  if and only if
(\ref{range}) holds.  If $p =0$, there exists a maximal $G$-invariant $K$-subalgebra
$R'$ of $R$ such that $R^{G} = R'^{G}$ and $(X', G)$ is stable for $X' = {\rm Spec} (R')$
(cf. \cite{Nak5-0}). 

In general if   $(X, G)$ is stable, 
${\mathcal Q}(R^G) = {\mathcal Q}(R)^G$ holds (cf. \cite{Po}), but 
the converse of this assertion does not hold as follows even if $G$ is an algebraic torus. 
\end{rem}

\begin{exam}\label{generalizedreflection} \rm 
The $K$-subalgebra   $B = K[X_{1}, X_{2}, X_{3}]$ of $A$ is invariant under 
the action of  $G= ({\mathbf{G_{m}}})^{2}$
 in Example \ref{stability} and ${\mathcal Q}(B^{G}) = {\mathcal Q}(B)^{G}$ holds. However the action $(Y,
G)$ is not   {\it stable} for $Y =  {\rm Spec} (B)$.  
\end{exam}

\begin{pr}\label{finitequotientcharactergroup}  Let $(S, G)$ be a regular action of 
  $G$  on  an  affine integral $K$-scheme $S = {\rm Spec} (A)$.   Let $L$ be a finite subgroup
of the centralizer $Z_{G}(G^{0})$ of $G^{0}$.   Then:
\begin{itemize}
\item[(i)] 
$\begin{aligned} & \sharp_{p'} (L\vert _{A}) \cdot {\frak X} (G^{0})^{A} \subseteq {\frak X} (G^{0})^{A^{L}}. 
\end{aligned}
$
\item[(ii)] $
\begin{aligned} &\sharp_{p'} (L\vert _{A}) \cdot {\frak X} (G^{0})_{A} \subseteq {\frak X} (G^{0})_{A^{L}}.
\end{aligned}$
\item[(iii)]  Both indices $[\langle{\frak X}(G^{0})^{A}\rangle : \langle{\frak X}(G^{0})^{A^{L}}\rangle]$ and $[{\frak X}(G^{0})_{A} : {\frak X}(G^{0})_{A^{L}}]$ are finite and  not divisible by
$p$ if $p >0$. 
\end{itemize}
\end{pr}

\proof  {\it (i):}~We suppose that $L$ acts {\it faithfully} on $A$. For any $\chi \in {\frak X} (G^{0})^{A}$, the non-zero $K$-subspace $A_{\chi}$  is invariant under the action of $L$, 
because  $L \subseteq Z_{G}(G^{0})$. Let $L_{p}$ 
denote a $p$-Sylow subgroup of $L$ if $p > 0$ and put  $L_{p} = \{1\}$ otherwise. Since 
the group $L_{p}$ is unipotent or trivial, we have a nonzero element  $x$ in
$A_{\chi}\cap A^{L_{p}}$. Expressing $$L = \bigsqcup_{i=1}^{u} \sigma_{i} L_{p}~\text{(a disjoint union)}$$ for $\sigma_{i} \in L$ $(1 \leqq i \leqq u = 
\sharp_{p'} (L)  )$,  put $$\tilde{x} = \prod_{i=1}^u {\sigma _{i}}(x) \in A_{u\chi}.$$
Then we easily see that $A^{L} \ni \tilde{x} \not= 0$, which implies $u\chi \in {\frak X}(G^{0})^{A^{L}}$
and {\it (i)}. 
 
 The assertion {\it (ii)} follows immediately from {\it (i)}. 
 Since the factor groups  $$ \langle{\frak X}(G^{0})^{A}\rangle / \langle{\frak X}(G^{0})^{A^{L}}\rangle ~\text{and} ~{\frak X}(G^{0})_A/ {\frak X}(G^{0})_{A^{L}}$$ are  finitely generated, by 
 {\it (i)} and  {\it (ii)} we see that the groups are  finite abelian groups with the exponents 
 which are  divisors of $\sharp_{p'} (L)$. Thus the assertion {\it (iii)} follows from this. 
 \qed

 \begin{lem}\label{commonquotient}  
  Let $(S_{i}, G)$ ($i =1, 2$) be  regular actions of  $G$ with
  an algebraic torus $G^{0}$ on affine
 integral $K$-schemes $S_{i} = {\rm Spec} (A_{i})$.  If ${\mathcal Q}(A_{1}) = {\mathcal Q}(A_{2}) $, then
 $$\langle {\frak X}(G^{0})^{A_{1}} \rangle = \langle {\frak X}(G^{0})^{A_{2}} \rangle .$$
 \end{lem}
 \proof Since the assertion is symmetry, we need only to show the inclusion
 $${\frak X}(G^{0})^{A_{2}}\subseteq \langle{\frak X}(G^{0})^{A_{1}}\rangle.$$
 Suppose that $(A_{2})_{\psi} \not= \{0\}$ for a character $\psi \in {\frak X}(G^{0})$ and
 let $b$ be a nonzero element of $(A_{2})_{\psi}$.  Let $(A_{1} : b)_{A_{1}}$ denote
 the ideal quotient $\{ a\in A_{1} \mid a\cdot b \in A_{1}\}$. As $b \in {\mathcal Q}(A_{1})
 = {\mathcal Q}(A_{2})$, 
we see that $(A_{1}: b)_{A_{1}}$ is a
 nonzero $G^{0}$-invariant ideal of $A_{1}$. This implies that 
 $$(A_{1})_{\chi}\supseteq ((A_{1} : b)_{A_{1}})_{\chi} \not=\{0\}$$
 for some $\chi \in {\frak X}(G^{0})$ and $$(A_{1})_{\psi + \chi}\supseteq
 b\cdot ((A_{1} : b)_{A_{1}})_{\chi} \not=\{0\}.$$ Thus 
 $$\psi = (\psi +\chi) - \chi \in \langle{\frak X}(G^{0})^{A_{1}}\rangle,$$ which
 shows ${\frak X}(G^{0})^{A_{2}}\subseteq \langle{\frak X}(G^{0})^{A_{1}}\rangle$. \qed
 
 Applying Lemma \ref{commonquotient} to Proposition \ref{finitequotientcharactergroup},
 we must have  
\begin{pr}\label{commonfinitegroupquotient}
Let $(S_{i}, G)$ ($i =1, 2$)  be  regular actions of  $G$ with
an algebraic torus $G^{0}$ on affine
 integral $K$-schemes $S_{i} = {\rm Spec} (A_{i})$. Let $L$ be a finite subgroup
of $Z_{G}(G^{0})$.   Suppose that ${\mathcal Q} (A_{1}^{L}) = 
{\mathcal Q}(A_{2})$. Then:
\begin{itemize}
\item[(i)] 
$\begin{aligned} &\sharp_{p'} (L\vert _{A_{1}}) \cdot \langle {\frak X} (G^{0})^{A_{1}}
\rangle \subseteq \langle {\frak X} (G^{0})^{A_{2}}\rangle 
.\end{aligned}$
\item[(ii)]  $\langle {\frak X} (G^{0})^{A_{1}}
\rangle$ contains  $\langle {\frak X} (G^{0})^{A_{2}}
\rangle$ as a subgroup whose index 
$$[\langle {\frak X} (G^{0})^{A_{1}}
\rangle : \langle {\frak X} (G^{0})^{A_{2}}
\rangle]$$ is not divisible by
$p$ if $p >0$. \qed
\end{itemize}

\end{pr} 

\begin{pr}\label{key}  Let $(X, G)$ be a regular action of  $G$
with an algebraic torus $G^{0}$ on an 
 affine integral   $K$-scheme $X= {\rm Spec} (R)$ such that $R$ is integrally closed.  Let  $L$ be  a finite subgroup of $Z_G(G^0)$.  Let $\mathfrak P$ be a $G^{0}$-invariant non-zero prime
ideal of $R$.  
Then $R^{L}$, $R^{{\mathcal I}_L({\mathfrak P})}$ and their reductions modulo by the ideals
induced by $\mathfrak P$ are naturally $G^{0}$-modules and the following properties hold:
\begin{itemize}
\item[(i)] $\sharp_{p'}({\mathcal D}_L({\mathfrak P})/{\mathcal I}_L({\mathfrak P}))\cdot \langle {\mathfrak X}(G^0)^{R^{{\mathcal I}_L({\mathfrak P})}/{\mathfrak P}^{{\mathcal I}_L({\mathfrak P})}}\rangle\subseteq
\langle{\mathfrak X}(G^0)^{{R^L}/{\mathfrak P}\cap R^L}\rangle$
\item[(ii)] $[\langle {\mathfrak X}(G^0)^{R^{{\mathcal I}_L({\mathfrak P})}/{\mathfrak P}^{{\mathcal I}_L({\mathfrak P})}}\rangle : \langle{\mathfrak X}(G^0)^{R^L/{\mathfrak P}\cap R^L}\rangle]$ is not divisible by $p$ if $p>0$. 
\item[(iii)] 
$\sharp_p\left(\langle {\mathfrak X}(G^0)^{R^{{\mathcal I}_L({\mathfrak P})}}\rangle/ \langle {\mathfrak X}(G^0)^{R^{{\mathcal I}_L({\mathfrak P})}/{\mathfrak P}^{{\mathcal I}_L({\mathfrak P})}}\rangle\right) = \sharp_p\left(\langle{\mathfrak X}(G^0)^{R^L}\rangle/\langle{\mathfrak X}(G^0)^{R^L/{\mathfrak P}\cap R^L}\rangle\right).$
\end{itemize}
\end{pr}
\proof  
The first assertion is obvious. 
Clearly  ${\mathcal D}_L({\mathfrak P})\vert _{R}$
is the decomposition group and ${\mathcal I}_L({\mathfrak P})\vert _{R}$ is the inertia group
at $\mathfrak P$ under the faithful action of $L\vert _{R}$ on $R$.   
The inclusions $$R^{L}  \subset R^{{\mathcal D}_L({\mathfrak P})} \subset R^{{\mathcal I}_L({\mathfrak P})}$$ induce the monomorphisms
 $$(R^{L}/{\mathfrak P}\cap R^{L}) \hookrightarrow
  (R^{{\mathcal D}_L({\mathfrak P})}/{\mathfrak P}\cap R^{{\mathcal D}_L({\mathfrak P})}) 
 \hookrightarrow \left(\left(R^{{\mathcal I}_L({\mathfrak P})}
 /{\mathfrak P}\cap R^{{\mathcal I}_L({\mathfrak P})}\right)^{{\mathcal D}_L({\mathfrak P})/{\mathcal I}_L({\mathfrak P})} \right).$$ 
 By the Galois theory $(R, L_{\vert R})$ of an integrally closed domain(cf. Chap. V of \cite{LR}), we see
 \begin{align*} & R^{{\mathcal D}_{L}({\mathfrak P})}_{{\mathfrak P}\cap  R^{{\mathcal D}_{L}({\mathfrak P})}}
 /({\mathfrak P}\cap  R^{{\mathcal D}_{L}({\mathfrak P})})R^{{\mathcal D}_{L}({\mathfrak P})}_{{{\mathfrak P}\cap  R^{{\mathcal D}_{L}({\mathfrak P})} }} \\ & \hskip 1.5cm =  \left(R^{{\mathcal I}_L({\mathfrak P})}_{{\mathfrak P}\cap  R^{{\mathcal I}_L({\mathfrak P})}}
 /({\mathfrak P}\cap  R^{{\mathcal I}_L({\mathfrak P})})R^{{\mathcal I}_L({\mathfrak P})}
 _{{\mathfrak P}\cap  R^{{\mathcal I}_L({\mathfrak P})}}\right)
 ^{{\mathcal D}_L({\mathfrak P})/{\mathcal I}_L({\mathfrak P})}\end{align*}
 and 
 $$R^{L}_{{\mathfrak P}\cap R^{L}}/({\mathfrak P}\cap R^{L}){R^{L}_{{\mathfrak P}\cap R^{L}}}=
 R^{{\mathcal D}_{L}({\mathfrak P})}_{{\mathfrak P}\cap  R^{{\mathcal D}_{L}({\mathfrak P})} }
 /({\mathfrak P}\cap  R^{{\mathcal D}_{L}({\mathfrak P})})R^{{\mathcal D}_{L}({\mathfrak P})}_{{{\mathfrak P}\cap  R^{{\mathcal D}_{L}({\mathfrak P})} }}. $$
These equalities imply 
 \begin{eqnarray*} {\mathcal Q}(R^{L}/{\mathfrak P}\cap R^{L}) & = &
 {\mathcal Q} (R^{{\mathcal D}_L({\mathfrak P})}/{\mathfrak P}\cap R^{{\mathcal D}_L({\mathfrak P})}) \\
 & = &
 {\mathcal Q}\left(R^{{\mathcal I}_L({\mathfrak P})}
 /{\mathfrak P}\cap R^{{\mathcal I}_L({\mathfrak P})}\right)^{{\mathcal D}_L({\mathfrak P})}\\
 & = &
{\mathcal Q}\left(\left(R^{{\mathcal I}_L({\mathfrak P})}
 /{\mathfrak P}\cap R^{{\mathcal I}_L({\mathfrak P})}\right)^{{\mathcal D}_L({\mathfrak P})}\right).\end{eqnarray*}
 Applying  Proposition \ref{commonfinitegroupquotient} to the replacements 
 $$A_{1}= R^{{\mathcal I}_L({\mathfrak P})}
 /{\mathfrak P}\cap R^{{\mathcal I}_L({\mathfrak P})}, ~A_{2}= R^{L}/{\mathfrak P} \cap R^{L}~ \text{and} ~L = {\mathcal D}_L({\mathfrak P}),$$ 
we obtain {\it (i)} and see that
$$\sharp_{p}\left(\langle {\mathfrak X}(G^0)^{R^{{\mathcal I}_L({\mathfrak P})}/{\mathfrak P}^{{\mathcal I}_L({\mathfrak P})}}\rangle / \langle{\mathfrak X}(G^0)^{R^L/{\mathfrak P}\cap R^L}\rangle]\right) =1,  $$
which shows {\it (ii)}.  On the other hand, applying 
Proposition \ref{commonfinitegroupquotient} to the replacement 
 $$A_{1}= R,~ A_{2}= R^{L}~\text{and}~ L = L,$$ we see also that $[\langle {\frak X} (G^{0})^{R}
\rangle : \langle {\frak X} (G^{0})^{R^{L}}
\rangle]$ and so 
$$[\langle {\mathfrak X}(G^0)^{R^{{\mathcal I}_L({\mathfrak P})}}\rangle : 
 \langle {\frak X} (G^{0})^{R^{L}}
\rangle]$$ are not divisible by
$p$ if $p >0$.

By the commutative diagram 
\[\xymatrix{
R^{{\mathcal I}_L({\mathfrak P})} \ar[r]^-{\text{can.}} & 
 R^{{\mathcal I}_L({\mathfrak P})}/{\mathfrak P}^{{\mathcal I}_L({\mathfrak P})} 
 \ar[r] & 0\\
R^{L} \ar@{^{(}->}[u] \ar[r]^-{\text{can.}}  & R^{L}/{\mathfrak P}\cap R^{L}
\ar@{^{(}->}[u]  \ar[r] & 0
}\]
with $G^{0}$-invariant exact rows, we have the descending chains
$$\left\{ \begin{aligned}
& \langle{\frak X} (G^{0})^{R^{{\mathcal I}_L({\mathfrak P})}}\rangle \supset
 \langle{\frak X} (G^{0})^{R^{{\mathcal I}_L({\mathfrak P})}/{\mathfrak P}^{{\mathcal I}_L({\mathfrak P})}}\rangle \supset \langle{\frak X} (G^{0})^{R^{L}/{\mathfrak P}\cap R^{L}}\rangle \\
 & \langle{\frak X} (G^{0})^{R^{{\mathcal I}_L({\mathfrak P})}}\rangle\supset
 \langle{\frak X} (G^{0})^{R^{L}}\rangle \supset \langle{\frak X} (G^{0})^{R^{L}/{\mathfrak P}\cap R^{L}}\rangle \end{aligned}\right.$$
 of subgroups. 
 Consequently the assertion {\it (iii)} follows from {\it (ii)} and the last conclusion of the former paragraph of the proof. \qed
 
In Sect.6 we use    Proposition \ref{key}   for an affine Krull 
$K$-scheme $X = {\rm Spec} (R)$ and any prime ideal  ${\mathfrak P} \in {\rm Ht}_{1}(R, R^{G})$.  In 
this case, ${\mathfrak P}$ is known to be  $G^{0}$-invariant (cf.  \cite{Magid}).

  \small \section{Finite Group Actions  on Toric Quotients}\normalsize

Throughout this section,  we suppose $G^{0}$
is an
 algebraic torus. 
For a subset $\Delta$ of $G$, let
$\mathrm{Unip} (\Delta)$ denote the set consisting of all elements in $\Delta$
which are unipotent in $G$.

\begin{pr}\label{lemmaunip} Let   $\Gamma$ be a subgroup 
$\langle {\rm Unip}(\Gamma) \cup \Gamma_{0}\rangle$  of $G$
for a finite subset $\Gamma_{0}$ of $G$ of finite orders.
 If $G = Z_{G}(G^{0})$, 
then  $\Gamma$ is finite.

\end{pr} 

\proof  Let $\rho : G \to GL(V)$ be a faithful finite dimensional rational representation 
of $G$.  The $K$-space $V$ is a direct sum $\bigoplus _{i=1} ^{n} V_{\chi_{i}}$  
of non-zero subspaces 
$V_{\chi_{i}}$ of $V$ of relative invariants (cf. {\bf 1.F.}) with respect to distinct  characters $\chi_{i}\in {\mathfrak X} (G^{0}) $. 
As $G \subseteq  Z_{G}(G^{0})$, each subspace $V_{\chi_{i}}$ ($1\leq i\leq n$) is 
invariant under the action of $G$ on $V$. Let $$\rho_{i} : G \to GL(V_{\chi_{i}})$$ denote
the representation defined in a natural way.  Since each element of $\rho_{i}(\mathrm{Unip} (G))$
is unipotent on $V_{\chi}$,  we have $$\mathrm{det}_{V_{\chi_{i}}}(\rho_{i}(\langle \mathrm{Unip} (G)\rangle)) =1$$ where $\mathrm{det}_{V_{\chi_{i}}}$ denotes the determinant on  $V_{\chi}$.  From this observation  and  the assumption on generators of $\Gamma$, 
we infer  that the image of the morphism 
$$\mathrm{det}_{V_{\chi_{i}}} : \rho_{i}(\Gamma)
{\longrightarrow} K^{\times}$$ is finite. For a nontrivial character
$\chi_{i}$, as $\rho_{i}(G^{0})$ is of rank one and the morphism $$\mathrm{det}_{V_{\chi_{i}}} : \rho_{i}(G^{0}) \to K^{\times}$$ is surjective, we see that 
$$\mathrm{Ker} (  \rho_{i}(G^{0}) \underset{\mathrm{det}_{V_{\chi_{i}}}}{\longrightarrow} K^{\times})$$ is a finite group. Consequently  the kernel and image of
morphism $$\rho_{i}(G^{0}{\cap }\Gamma)
\underset{\mathrm{det}_{V_{\chi_{i}}}}{\longrightarrow} K^{\times}$$
are finite,  which imply $\rho_{i}(G^{0}{\cap }\Gamma)$ is a finite for a nontrivial $\chi_{i}$. 
As   $\rho_{i} (\Gamma)/\rho_{i}(G^{0}\cap \Gamma)$ is  the epimorphic image
of the finite group $\Gamma/G^{0}\cap \Gamma$, it  is finite if $\chi _{i}$ is nontrivial. On 
the other hand, if  $\chi_{i}$ is the trivial character, $\rho_{i}(G)$ is clearly finite. 
Because 
 $\Gamma$ can be embedded in  the product  $$ \prod_{i=1}^{n} \rho_{i}(\Gamma)$$
of finite groups,
we conclude that $\Gamma$ 
is also finite. 
 \qed

\paragraph{\bf Notation} 
Hereafter in this section, let  $(X, G)$ denote a regular effective action of 
$G$ on an affine Krull $K$-scheme $X = {\rm Spec} (R)$. 
Clearly $G\vert_{R^{G^{0}}}$ is a finite group. Hence, let  ${\mathcal  I}_G({\mathfrak q})_{[p]}$ denote the inverse image of  a  $p$-Sylow subgroup of
${\mathcal  I}_G({\mathfrak q})\vert_{R^{G^0}}$ 
under the canonical homomorphism 
$${\mathcal  I}_G({\mathfrak q})\to {\mathcal  I}_G({\mathfrak q})\vert_{R^{G^0}}$$ for any ${\mathfrak q}\in {\rm Ht}_1(R^{G^0})$ 
if $p ~( = {\rm char} (K)) > 0$.   Otherwise put  $${\mathcal  I}_G({\mathfrak q})_{[0]} = 
{\mathrm{Ker}} ({\mathcal I}_G({\mathfrak q})
\to   {\mathcal I}_G({\mathfrak q})\vert_{R^{G^0}}).$$
Both groups  ${\mathcal I}_G({\mathfrak q})$ and  ${\mathcal  I}_G({\mathfrak q})_{[p]}$  contain $G^{0}$ and are  closed in $G$. 

\begin{lem}\label{lemma2-6}  Let ${\mathfrak q}$ be any prime ideal
in ${\rm Ht}_1(R^{G^0})$.   Then ${\mathcal I}_G({\mathfrak q})_{[p]}$ is normal in ${\mathcal I}_G({\mathfrak q})$
and $ {\mathcal I}_G({\mathfrak q})/{\mathcal I}_G({\mathfrak q})_{[p]}$ is a subgroup of 
of ${\rm U}(R^{G^0}_{\mathfrak q}/{\mathfrak q}R^{G^0}_{\mathfrak q})$. 
\end{lem} 

 \proof Recalling   ${\mathcal  I}_G({\mathfrak q}) = 
\{ \sigma \in G \mid (\sigma - 1)(R^{G^0}) \subseteq {\mathfrak q}\}$,  we see
$${\mathcal  I}_G({\mathfrak q})\vert_{R^{G^0}} = {\mathcal I}_{G\vert _{R^{G^{0}}}}({\mathfrak q})
= {\mathcal I}_{G\vert _{R^{G^0}}}(R^{G}_{R^{G}\cap{\mathfrak q}}\otimes _{R^{G}}{\mathfrak q})$$ which is the
inertia group of $R^{G}_{R^{G}\cap{\mathfrak q}}\otimes _{R^{G}}{\mathfrak q}$ for  the finite Galois group action 
$$(R^{G}_{R^{G}\cap{\mathfrak q}}\otimes _{R^{G}}R^{G^0}, G\vert _{R^{G^0}})$$ on the semilocal Dedekind domain.  
We now examine  and apply the proof of 
 Theorem 25 in  \S 10 of Chap. V of  \cite{zariski} to our case. Then the $p$-Sylow
subgroup $H_{1}$ of ${\mathcal I}_G({\mathfrak q})\vert_{R^{G^0}}={\mathcal I}_{G\vert _{R^{G^0}}}(R^{G}_{R^{G}\cap{\mathfrak q}}\otimes _{R^{G}}{\mathfrak q})$
is normal  if $p >0$. Regarding $H_{1}$ is trivial in the case where
$p=0$ we always have $$
 ({\mathcal I}_G({\mathfrak q})\vert_{R^{G^0}})/H_{1}\hookrightarrow {\rm U}(R^{G^0}_{\mathfrak q}/{\mathfrak q}R^{G^0}_{\mathfrak q})$$
as groups for any characteristic 
$p$. Thus the assertion follows from this observation.   \qed 

\begin{pr}\label{pr2-7}  Suppose that  $G = Z_G(G^0)$. 
 For any ${\mathfrak q}\in {\rm Ht}_1(R^{G^0})$, we have:
 \begin{itemize}
 \item[(i)] $\langle {\rm Unip} ({\mathcal  I}_G({\mathfrak q}))\rangle $ is a finite group satisfying
$\langle{\rm Unip} ({\mathcal  I}_G({\mathfrak q}))\rangle \vert_{R^{G^0}} =
   {\mathcal  I}_G({\mathfrak q})_{[p]} \vert_{R^{G^0}}$ and 
   $$R^{\langle {\rm Unip} ({\mathcal  I}_G({\mathfrak q}))\rangle \cdot G^0} =
   R^{{\mathcal  I}_G({\mathfrak q})_{[p]}}.$$ 
   \item[(ii)]  There exists
   a diagonalizable closed subgroup,  denoted by ${\mathcal D}_{G, {\mathfrak q}}$,  of ${\mathcal I}_G({\mathfrak q})$
   such that $D _{G, {\mathfrak q}}\supseteq G^0$, ${\mathcal D}_{G, {\mathfrak q}}/G^0$ is cyclic, and $$R^ {\langle {\rm Unip} ({\mathcal  I}_G({\mathfrak q}))\rangle \cdot {\mathcal D}_{G, {\mathfrak q}}}= R^{{\mathcal I}_G({\mathfrak q})}.$$
   \item[(iii)]   $R^{G}_{{\mathfrak  q}\cap R^{G}} \to R^{{\mathcal I}_G({\mathfrak q})}
   _{{\mathfrak  q}\cap R^{{\mathcal I}_G({\mathfrak q})}}$ is \'etale  and the ramification index $$\mathrm{e}( {\mathfrak q}\cap R^{\langle {\rm Unip} ({\mathcal  I}_G({\mathfrak q}))\rangle \cdot G^0}, {\mathfrak q}\cap R^G)$$ is not divisible by $p$ if $p >0$. 
   
   \end{itemize} 
\end{pr} 

\proof  {\it (i):} ~By Proposition \ref{lemmaunip} we easily see  that  the group $ \langle {\rm Unip} (G)\rangle$ is of finite order. 
Let $\sigma \in  {\mathcal  I}_G({\mathfrak q})_{[p]}$ and  express $\sigma = \sigma_u \cdot \sigma_s$
for a unipotent (resp. semisimple) element $\sigma_u$ (resp. $\sigma_s$)  in ${\mathcal  I}_G({\mathfrak q})_{[p]}$. 
As $\sigma \vert_{R^{G^0}}$ is unipotent,  
we see $$\sigma \vert_{R^{G^0}}=
\sigma_u \vert_{R^{G^0}}\in {\rm Unip} ({\mathcal  I}_G({\mathfrak q})) \vert_{R^{G^0}}.$$  On the other hand, the
inclusion $${\rm Unip} ({\mathcal  I}_G({\mathfrak q})) \vert_{R^{G^0}} \subseteq
   {\mathcal  I}_G({\mathfrak q})_{[p]}\vert_{R^{G^0}}$$ follows, since 
    the number of $p$-Sylow subgroups of  ${\mathcal I}_G({\mathfrak q})\vert_{R^{G^0}} $ is one  (cf. Lemma \ref{lemma2-6}) if $p >0$. Consequently 
     $$\langle{\rm Unip} ({\mathcal  I}_G({\mathfrak q}))\rangle \vert_{R^{G^0}} =
   {\mathcal  I}_G({\mathfrak q})_{[p]} \vert_{R^{G^0}}$$ holds.  From this we immediately derive
   the second equality of {\it (i)}. 
   
   {\it (ii):}~ According to Lemma \ref{lemma2-6}, choose  $\tau \in {\mathcal I}_G({\mathfrak q})$ in such a way that $\tau$  generates the factor group 
   $$I_G({\mathfrak q})\vert_{R^{G^0}}/({\mathcal  I}_G({\mathfrak q})_{[p]})\vert_{R^{G^0}}$$ on $R^{G^{0}}$ principally. 
   Expressing $\tau =\tau_u \cdot \tau_s$ with its  unipotent part  $\tau_u$ and semisimple one  $\tau_s$,
   we may assume $\tau$ is semisimple.   Define ${\mathcal D}_{G, {\mathfrak q}}$ to be the subgroup $\langle G^0, \tau\rangle$.
    We must have 
   $$R^ {\langle {\rm Unip} ({\mathcal  I}_G({\mathfrak q}))\rangle \cdot {\mathcal D}_{G, {\mathfrak q}}}= 
   (R^{G^0})^{{\langle {\rm Unip} ({\mathcal  I}_G({\mathfrak q}))\rangle \cdot {\mathcal D}_{G, {\mathfrak q}}}/G^0} = 
   (R^{G^0})^{{\mathcal I}_G({\mathfrak q})/G^0} = R^{{\mathcal I}_G({\mathfrak q})}$$
   since $\langle {\rm Unip} ({\mathcal  I}_G({\mathfrak q}))\rangle  \vert_{R^{G^0}} =
   {\mathcal  I}_G({\mathfrak q})_{[p]}\vert_{R^{G^0}}$.
   
   {\it (iii):}  By the Galois theory on
   $R^{G^{0}}$ with the finite group  $G/G^{0}$ action and Corollary \ref{etale2}, we see that
   $$R^{G} \to R^{{\mathcal I}_G({\mathfrak q})}$$ is \'etale at 
   ${\mathfrak  q}\cap R^{{\mathcal I}_G({\mathfrak q})}$. Suppose that $p > 0$ and consider  
   $$(R^{\langle {\rm Unip} ({\mathcal  I}_G({\mathfrak q}))\rangle \cdot G^0})
   _{{\mathfrak q}\cap R^{\langle {\rm Unip} ({\mathcal  I}_G({\mathfrak q}))\rangle \cdot G^0}}$$ acted by ${\mathcal D}_{G, {\mathfrak q}}/G^{0}$ in  {\it (ii)} of this proposition. 
  Applying  the ramification theory 
   of Dedekind domains (e.g., Chap. V of \cite{zariski})  to this discrete valuation ring, 
   $$\mathrm{e}( {\mathfrak q}\cap R^{\langle {\rm Unip} ({\mathcal  I}_G({\mathfrak q}))\rangle \cdot G^0}, {\mathfrak q}\cap R^{{\mathcal  I}_G({\mathfrak q})})$$
  is a divisor of  the order of the group 
  $$({\mathcal D}_{G, {\mathfrak q}}/G^{0})\vert _{R^{\langle {\rm Unip} ({\mathcal  I}_G({\mathfrak q}))\rangle \cdot G^0}},$$ which is not divisible by $p$.  Thus the assertion {\it (iii)} follows from this. \qed

\section{Towards Krull Schemes}
We now generalize the auxiliary results in Sec. 3 of \cite{Nak2} for affine normal varieties
to  affine Krull $K$-schemes.  Stability of actions 
$(X, G)$ with tori $G^{0}$ can be replaced by ${\mathcal Q}({\mathcal O}(X)^{G})
= {\mathcal Q}({\mathcal O}(X))^{G}$. 

Moreover we need in the last section  the following theorem which is the main result of \cite{Nak5} and
implies that  {\it any inertia group ${\mathcal I} _{G}({\mathfrak Q})$ at ${\mathfrak Q} \in \mathrm{Ht}_{1}({\mathcal O}(X), {\mathcal O}(X)^G)$ 
is finite on $X$ for a regular action$(X, G)$  of a reductive algebraic group $G$  on 
an arbitrary affine Krull $K$-scheme $X$. }
Pseudo-reflection groups
${\mathfrak R}({\mathcal O}(X), G)$ of actions   are referred to the paragraph {\bf 2.D.} of Sect. 2.

\begin{theorem}[cf. \cite{Nak5}]\label{mae}  In general, the following conditions are equivalent:
\begin{itemize}
\item[(i)]  $G$ is (geometrically) reductive.
\item[(ii)] For any regular action $(X,G)$ of $G$ on any 
affine Krull $K$-scheme $X$, the image of the pseudo-reflection group
${\mathfrak R}({\mathcal O}(X), G)$ to ${\rm Aut} X$ (i.e.,  ${\mathfrak R}({\mathcal O}(X), G)\vert_{{\mathcal O}(X)}$) is finite. \qed
\end{itemize}
\end{theorem}

The next result, which is a slight generalization of Proposition 1.8 of \cite{Nak2},  can be shown as in the proof in  \cite{Nak2}
without the assumption that $G^{0}$ is {\it linearly} reductive  and
$X$ is noetherian by the use of Theorem \ref{mae} and
Corollary \ref{etale2}.

\begin{pr}[cf.  \cite{Nak2}]\label{reflection} Suppose that $G$ is a reductive affine algebraic 
and   that $X$ is an affine Krull $K$-scheme with a 
regular action $(X, G)$ of $G$. 
Set $Y = X\dslash {\mathfrak R}({\mathcal O}(X), G)$. Then :
\begin{itemize}
\item[(i)]  ${\mathfrak R} ({\mathcal O}(Y) , G) = {\mathfrak R} ({\mathcal O}(X), G)$.
\item[(ii)] ${\mathcal I}_{G}({\mathfrak q})\vert_{Y} = \{1\}$
for ${\mathfrak q} \in {\rm Ht}_{1}({\mathcal O}(Y), {\mathcal O}(Y)^G)$. \qed
\end{itemize}
\end{pr}

\begin{pr}\label{keylemma} 
 Let $D$ be a diagonalizable affine algebraic group 
and $(X, D)$  a regular action of $D$ on an  affine Krull $K$-scheme   $X = {\rm Spec} (R)$  such that
${\mathcal Q}(R^{D}) = {\mathcal Q}(R)^{D}$. 
Let $\mathfrak Q$ be a prime ideal in $\mathrm{Ht}_{1} (R, R^D)$.  If ${\mathfrak Q}$ is invariant under the action of $D$,  then $\langle {\mathfrak X}(D)^{R/{\mathfrak Q}}\rangle $ 
is a subgroup of  $\langle {\mathfrak X}(D)^{R}\rangle $
and the  quotient group 
 $$\langle {\mathfrak X}(D)^{R}\rangle/
 \langle {\mathfrak X}(D)^{R/{\mathfrak Q}}\rangle $$
 is a cyclic group of   order  equal to the ramification index
$\mathrm{e} ({\mathfrak Q}, {\mathfrak Q}\cap R^{D})$. 
\end{pr}
\proof  Since ${\mathfrak Q}$ is a rational $D$-submodule
and $D$ is linearly reductive, for each $\chi \in 
{\mathfrak X}(D)^{R/{\mathfrak Q}}$ we have $R_{\chi} \not= \{0\}$, which implies
 the first assertion and the equality 
\begin{equation} {\mathfrak X}(D)^{R} = {\mathfrak X}(D)^{{\mathfrak Q}}
\cup {\mathfrak X}(D)^{R/{\mathfrak Q}}. \label{characterdecomposition}\end{equation}

Let $f$ be
an element  of ${\mathfrak Q}$ to satisfy $f\cdot R_{\mathfrak Q} = {\mathfrak Q}R_{\mathfrak Q}$. As there is  a finite dimensional rational $D$-submodule $U$
of ${\mathfrak Q}$  containing $f$,  exchanging $f$, we can choose the element $f$ 
 in such a way that $f \in R_{\psi}$ for some $\psi \in {\mathfrak X}(D)^{R}$. Let
 $\mu$ be a character in ${\mathfrak X}(D)^{\mathfrak Q}$.  Then 
 we have a nonzero element $y\in {\mathfrak Q}$ such that $y \in R_{\mu}$.  There
 are 
 nonzero elements $g, h \in R \backslash {\mathfrak Q}$ satisfying 
\begin{equation} y = f^{{\rm v}_{R, {\mathfrak Q}}(y)}\cdot  \frac{g}{h}.
\label{eq:y}\end{equation}
 Express the finite sums $$g= \sum_{\xi \in {\mathfrak X}(D)} g_{\xi} ~~(g_{\xi}\in R_{\xi}) ~
 \text{and} ~
h= \sum_{\xi \in {\mathfrak X}(D)} h_{\xi} ~~(h_{\xi}\in R_{\xi}).$$ Then
 $h_{\xi_{1}} \not\in {\mathfrak Q}$ for some $\xi_{1} \in {\mathfrak X}(D)$. 
 Clearly by (\ref{eq:y}) we have  $$y \cdot h_{\xi_{1}} =  f^{v_{R, {\mathfrak Q}}(y)}\cdot g_{\xi_{2}}$$
 for some $\xi _{2}\in {\mathfrak X}(D)$. 
 Comparing the values of ${\rm v}_{R, {\mathfrak Q}}$ of the
 above equality, we must have 
 $g_{\xi_{2}}\not\in{\mathfrak Q}$, which implies 
 $$\mu \equiv {v_{R, {\mathfrak Q}}(y)} \psi \mod \langle {\mathfrak X}(D)^{R/{\mathfrak Q}}\rangle$$
 because $\xi_{1}, \xi_{2} \in  {\mathfrak X}(D)^{R/{\mathfrak Q}}$. 
 Consequently by (\ref{characterdecomposition}) we see $$\langle {\mathfrak X}(D)^{R}\rangle/
 \langle {\mathfrak X}(D)^{R/{\mathfrak Q}}\rangle $$ is a cyclic
 group generated by $\psi + \langle {\mathfrak X}(D)^{R/{\mathfrak Q}}\rangle$. 
 
 By Lemma  \ref{quotient},  we have $R_{\mathfrak Q} \cap {\mathcal Q}(R)^{D}
 =R^{D}_{{\mathfrak Q}\cap R^{D}}$ and 
 $$f R_{\mathfrak D} \cap {\mathcal Q}(R)^{D}= {\mathfrak Q}R_{\mathfrak Q} \cap {\mathcal Q}(R)^{D}
 =({{\mathfrak Q}\cap R^{D}})R^{D}_{{\mathfrak Q}\cap R^{D}}.$$
 From  ${\rm v}_{R, {\mathfrak Q}}({\mathfrak Q}\cap R^{D}) = 
 \mathrm{e} ({\mathfrak Q}, {\mathfrak Q}\cap R^{D})$, we infer that 
 \begin{equation*} \mathrm{e} ({\mathfrak Q}, {\mathfrak Q}\cap R^{D})
 = \mathrm{min} \left\{ k \in \mbox{\boldmath $N$}~ \left\vert~
 \exists g, h \in R\backslash {\mathfrak Q} ~\text{such that} ~ 
 f^{k}\cdot \frac{g}{h} \in  {\mathcal Q}(R)^{D}
 \right. \right\}\end{equation*}
 As in the preceding paragraph, this can be   replaced by 
 \begin{equation}\label{e} \mathrm{e} ({\mathfrak Q}, {\mathfrak Q}\cap R^{D}) =
  \mathrm{min} \left\{ k \in \mbox{\boldmath $N$}~\left\vert~ \begin{array}{c} 
 \exists ~\text{relative invariants} ~ g, h \in R\backslash {\mathfrak Q} ~ \text{of}~ D   \\ \text{such that} ~
  f^{k}\cdot \displaystyle\frac{g}{h} \in  {\mathcal Q}(R)^{D}\end{array}\right.
 \right\}\end{equation}
 (see {\bf 1.F.} in Sect. 1 for {\it relative invariants} of $D$). 
 For a rational character $\chi \in \langle {\mathfrak X}(D)^{R}\rangle$, 
 $\chi \in  \langle {\mathfrak X}(D)^{R/{\mathfrak Q}}\rangle$ if and only if
 \begin{align*} \chi  & = \sum_{\exists \{\xi_{i}\}\subset  {\mathfrak X}(D)^{R/{\mathfrak Q}}}  \xi_{i} - 
  \sum_{\exists \{\nu_{j}\} \subset {\mathfrak X}(D)^{R/{\mathfrak Q}}} \nu_{j} ~~\text{(finite sum)}\\
  & = \xi - \nu ~~~(\exists \xi, \nu \in {\mathfrak X}(D)^{R/{\mathfrak Q}}), 
  \end{align*}
  because the integral domain $R/{\mathfrak Q}$ implies that ${\mathfrak X}(D)^{R/{\mathfrak Q}}$ is a semi-group. 
 Thus the right hand side of (\ref{e}) can be identified with 
 $$\mathrm{min} \left\{  k \in \mbox{\boldmath $N$}~\left\vert~
 k\mu \in   \langle {\mathfrak X}(D)^{R/{\mathfrak Q}}\rangle \right. \right\},$$
 which is equal  to the order of the group 
 $\langle {\mathfrak X}(D)^{R}\rangle / \langle {\mathfrak X}(D)^{R/{\mathfrak Q}}\rangle$.
 \qed 
 
 If a regular action  $(X, G)$ is effective (cf. Sect. 1),   
 \begin{equation} {\rm Ker} (G \to {\rm Aut}~X)
\subseteq Z_G(G^0)\label{ineffectivekernel}\end{equation} follows from its finiteness.
A variety $W$ is said to be {\it conical} if $X$ is affine and ${\mathcal O}(W)$ is
a graded $K$-algebra equipped with a positive gradation
${\mathcal O}(W) = \oplus_{i\geq 0}{\mathcal O}(W)_{i}$  such that ${\mathcal O}(W)_{0} = K$.  In this case 
a regular action $(W, G)$  is said to be conical, if it is regular and 
the induced action of $G$ on ${\mathcal O} (W)$ preserves the
gradation of ${\mathcal O}(W)$. 
The following result is the restatement of  Proposition 2.2 of \cite{Nak2}
for $W = \tilde{V}\dslash L$ in this proposition.  
 
\begin{pr}[cf. \cite{Nak2},  Propositions 2.2]\label{Nak22}  Suppose that $G^{0}$
is an algebraic torus.  Then 
there exists an effective conical stable action  $(W, G)$ 
of on  a conical normal variety $W$ 
such that $(W\dslash G^0, G)$ can be identified with  a trivial action
on the affine line $\mbox{\boldmath $A$}^1$.  \qed
\end{pr}

\begin{pr}[cf. \cite{Nak2}, Proposition 2.3]\label{Nak23}  Suppose that 
$G^{0}$ is  an algebraic torus and $G \not= Z_{G}(G^{0})$.
  If   $p = 0$ or $\sharp (G/Z_{G}(G^{0}))$ is relatively prime to $p$ for $p >0$, then 
 there exists a closed subgroup $H$ of $G$ containing $Z_{G}(G^{0})$ and an
effective  regular action $(X, H)$ 
of $H$ on  an affine normal variety $X$ such that 
 \begin{equation*} {\rm e}_{p'}({\mathfrak P}, {\mathfrak P}\cap {\mathcal O}(X)^H)
\not=
\sharp_{p'} ( {\mathcal I}_{H}({\mathfrak P})\vert _{X}).
\end{equation*}
\end{pr}

\proof  We choose a closed subgroup  $H$ of $G$ in such a way that 
$H/Z_{G}(G^{0})$ is non-trivial cyclic subgroup which is regarded as 
a subgroup of  the multiplicative group $\mathbf{G_{m}}$ over $K$.  Consider   the canonical regular action 
$(Y, \mathbf{G_{m}})$ with $Y = {\rm Spec} (K[X_{1}])$ 
on the affine line $\mbox{\boldmath $A$}^1$ such as $t(X_{1})= tX_{1} ~(t\in \mathbf{G_{m}})$. 
By this  we define
a regular action $(Y, H/Z_{G}(G^{0}))$.  Applying Proposition 2.3 to $(Y, H/Z_{G}(G^{0}))$, 
we obtain an effective  regular action $(X, H)$ on an affine normal variety $X$
as desired.  \qed 

Using   Proposition \ref{keylemma} and Proposition \ref{reflection}
instead of Lemma 3.1 (1)  and Proposition 1.8 of \cite{Nak2} respectively, 
we can similarly   show the following  result which is a slight generalization of 
  Proposition 3.2 of \cite{Nak2}.  

\begin{pr}\label{mainofpreviouspaper} Suppose that $G^{0}$ is an algebraic torus and $G = Z_{G}(G^{0})$. 
Let $(X, G)$ be a regular action of $G$ on an affine Krull $K$-scheme 
$X = {\rm Spec} (R)$ such 
that ${\mathcal Q}(R)^{G} = {\mathcal Q}(R^G)$. Then ${\rm e} ({\mathfrak Q}, {\mathfrak Q}\cap R^{G} )$ is equal to  
 \begin{align*}
\sharp_{{p'}}  ({\mathcal I}_{G}({\mathfrak Q})\vert _{R})\cdot 
\sharp_{p}\left(\langle{\mathfrak X}(G^{0})^{R}\rangle/\langle{\mathfrak X}(G^{0})^{R/{\mathfrak Q}}\rangle\right)
 \cdot {\rm e}_{p} ({\mathfrak Q}\cap R^{G^{0}},
 {\mathfrak Q}\cap R^{{\mathcal I}_{G}({\mathfrak Q}\cap G^{0})}) %
\end{align*} 
for any ${\mathfrak Q} \in \mathrm{Ht}_{1} (R, R^G)$.  \text{\qed}

\end{pr}

\section{Main Theorem}
As a matter of convenience, we restate Theorem 1.1 as follows: 
\begin{theorem}\label{mainresult}
Suppose that $G^0$ is an algebraic torus.  Then the following conditions are equivalent:
\begin{itemize}
\item[(i)]  $G = Z_G(G^0)$
\item[(ii)] For an arbitrary   closed subgroup $H$ of   $G$ containing $Z_G(G^0)$, 
the following 
conditions hold  for any 
effective regular action   $(X, H)$
on an arbitrary affine Krull $K$-scheme $X = {\rm Spec} (R)$ such that ${\mathcal Q}(R^{G^{0}}) = {\mathcal Q}(R)^{G^{0}}$:   for  any ${\mathfrak P}\in
{\rm Ht}_1(R, R^H)$, 
$${\rm e}({\mathfrak P}, {\mathfrak P}\cap R^H)={\rm e}({\mathfrak P}, {\mathfrak P}\cap R^{{\mathcal I}_H({\mathfrak P})}) \\
\cdot \sharp_p \left(
\langle {\mathfrak X}(G^0)^{R^{{\mathcal I}_{H}({\mathfrak P})}}\rangle /\langle{\mathfrak X}(G^0)^{R^{{\mathcal I}_{H}({\mathfrak P})}/{\mathfrak P}\cap R^{{\mathcal I}_{H}({\mathfrak P})}}\rangle\right).
$$
\item[(iii)] For an arbitrary   closed subgroup $H$ of   $G$ containing $Z_G(G^0)$, 
the 
conditions in {\it (ii)} hold  for any 
effective stable  regular action   $(X, H)$
on an arbitrary affine normal variety   $X = {\rm Spec} (R)$.

\end{itemize}
\end{theorem}

\proof  Since   $(X, H)$ in {\it (iii)} is a
special  case mentioned  in  {\it (ii)} (cf. Remark \ref{stabilitycriteria}), the implication {\it (ii)} $\Rightarrow${\it (iii)} is obvious. 
Note that  $Z_{H}(G^{0}) = Z_{G} (G^{0})$ for  a subgroup $H$ of $G$
containing $Z_{G}(G^{0})$ and that in {\it (ii)}  the equivalence 
$${\mathcal Q}({\mathcal O}(X)^{G^{0}}) = {\mathcal Q}({\mathcal O}(X))^{G^{0}}
\Longleftrightarrow  {\mathcal Q}({\mathcal O}(X)^{H}) = {\mathcal Q}({\mathcal O}(X))^{H}.$$ follows from Lemma \ref{qf}. Moreover in {\it (ii)}
note that   ${\mathcal I}_H({\mathfrak P})$ is a finite group.   

{\it (iii)}$\Rightarrow${\it (i)}:  Assume that this implication is false. Suppose
$G$ is a counter-example for  {\it (iii)}$\Rightarrow${\it (i)}, {\it i.e.}, the condition {\it (iii)} holds for $G$ but $G \not= Z_{G}(G^{0})$. 
By {\it (iii)} and (\ref{eq1}) we always 
have \begin{equation} {\rm e}_{p'}({\mathfrak P}, {\mathfrak P}\cap {\mathcal O}(X)^H)
= {\rm e}_{p'}({\mathfrak P}, {\mathfrak P}\cap R^{{\mathcal I}_H({\mathfrak P})}) =
\sharp_{p'} ( {\mathcal I}_{H}({\mathfrak P})\vert _{X})\label{false}
\end{equation}
under the circumstances as in {\it (ii)}. The second equality follows from (\ref{eq1}),  since
${\mathcal I}_{H}({\mathfrak P})$ is finite on $X$.
 From Proposition \ref{Nak23} we see 
 $p > 0$ and the  set  $G\backslash Z_{G}(G^0)$
 contains an element $\tau$ such that  order of  the coset $\tau Z_{G}(G^{0})$ in
 the group $G/Z_{G}(G^{0})$ is equal to $p$.   
 
 Let $$\phi: C_1 \to C_2$$ be a  Galois cover of smooth
 projective curves defined over $K$ which is wildly ramified at a closed  point (e.g., \cite{Pries2}).  Exchanging $C_1$ with
 its quotient by the normal subgroup of the Galois group of index $p$, we may suppose that the Galois group of $\phi$ is equal  to 
 $$\Gamma \cong \mbox{\boldmath $Z$}/p\mbox{\boldmath $Z$}.$$ Then there exists a  closed point  $P\in C_1$
 fixing by $\Gamma$, because the inertia group ${\mathcal I}_{\Gamma}(P)$ of $P$ under the action of $\Gamma$ must be equal to
 $\Gamma$.  Let $Y$ be an affine open set containing the closed point $\phi (P)$. Let $H$ be the closed subgroup of $G$
 generated by $Z_G(G^0)$ and $\tau$ whose quotient group
 $$H/Z_G(G^0) \cong \Gamma$$ acting regularly on the affine normal variety $\phi^{-1}(Y)$.
 So $\phi^{-1}(Y)\dslash H = Y$. 
 Let ${\mathfrak P}_{Y}$ be the maximal ideal  of ${\mathcal O}(\phi^{-1}(Y))$ associated with
 the point $P \in \phi^{-1}(Y)$.  Since $\sharp (H\vert_{\phi^{-1}(Y)}) = \sharp ( {\mathcal I}_{H}({\mathfrak P}_{Y})\vert_{\phi^{-1}(Y)}
) = \sharp( \Gamma ) = p$ and ${\mathcal O}(\phi^{-1}(Y))/{\mathfrak P}_{Y}
 = {\mathcal O}(Y)/({\mathfrak P}_{Y}\cap {\mathcal O}(Y)) = K$, by ramification theory
 (e.g., Chap. V of \cite{zariski}) we see that ${\rm e}({\mathfrak P}_{Y}, {\mathfrak P}_{Y}\cap {\mathcal O}(Y)) =p$.
 From Proposition \ref{Nak22} there exists a  conical effective action of $H$ 
 on a conical normal variety $W$ 
such that $(W\dslash G^0, H)$ can be identified with  a trivial action
on the affine line $\mbox{\boldmath $A$}^1$. 
Consider the affine normal variety 
 $$X := \phi^{-1}(Y) \times_{K} W$$ 
 with ${\mathcal O} (X) = {\mathcal O}(\phi^{-1}(Y)) \otimes_{K} {\mathcal O}(W)$ on which $H$  acts diagonally. As $H^{0} = G^{0}$ 
  acts on $\phi^{-1}(Y)$ trivially, the regular action $(X, H)$ is stable.  Put $${\mathfrak P} := {\mathfrak P}_{Y}\otimes_K{\mathcal O}(W) \in 
 {\rm Ht}_{1}({\mathcal O}(X)). $$
 Since 
 $${\mathcal O}(X)/{\mathfrak P} \cong {\mathcal O}(W)$$  is
 ${\mathcal D}_H({\mathfrak P})$-invariant (cf. {\bf 2.B.} of Sect. 2),  from the
 definition of  inertia groups and (\ref{ineffectivekernel}) we must have \begin{equation}\label{eqn6-1} {\mathcal I}_H({\mathfrak P}) \subseteq 
 {\rm Ker} (H \to {\rm Aut}~W) \subseteq Z_H(G^0).\end{equation}
 As the action of $Z_H(G^0)$ on $\phi^{-1}(Y)$ is trivial,  by (\ref{eqn6-1}) we see 
 \begin{equation}\label{eqn6-2}{\mathcal I}_H({\mathfrak P})\vert_X= \{1\}.\end{equation}
  On the other hand the triviality of $Z_H(G^0)$ on $\phi^{-1}(Y)$ and 
  $H$ on $W\dslash G^{0} \cong \mbox{\boldmath $A$}^{1}$ implies the
  commutative diagram with  vertical isomorphisms and horizontal quotient morphisms
 \[\xymatrix{
X\dslash G^0 \ar[d]^{\cong} \ar[r]^-{\text{can.}}& 
X\dslash Z_{H}(G^0)  \ar[d]^{\cong} \ar[r]^-{\text{can.}}& 
 \ar[d]^{\cong} X\dslash H
 \\
\phi^{-1}(Y) \times_{K} \mbox{\boldmath $A$}^1 \ar[r]^-{=}  & 
\phi^{-1}(Y) \times_{K} \mbox{\boldmath $A$}^1\ar[r]^-{\text{can.}}  &
(\phi^{-1}(Y)\dslash H )\times_{K}  \mbox{\boldmath $A$}^1= Y\times_{K} \mbox{\boldmath $A$}^1~.
} \]
 Since ${\mathfrak P}\cap {\mathcal O}(X) ^{G^{0}} = {\mathfrak P}_{Y}\otimes_{K} {\mathcal O}(W)^{G^{0}}$  and 
 ${\mathfrak P}\cap {\mathcal O}(X) ^{H} = ({\mathfrak P}_{Y}\cap {\mathcal O}(Y))\otimes_{K} {\mathcal O}(W)^{G^{0}}$,  we see that
$$ \rm{ht} ~({\mathfrak P}\cap {\mathcal O}(X) ^{H} ) = \rm{ht} ~({\mathfrak P}\cap {\mathcal O}(X) ^{G^{0}}) =1$$ and 
 \begin{equation}\label{eqn6-3} {\rm e}({\mathfrak P}\cap {\mathcal O}(X)^{G^0}, {\mathfrak P}\cap {\mathcal O}(X)^{H}) = p.\end{equation}
 Put $ M ={\mathcal O}(X)^{{\mathcal I}_{H}({\mathfrak P})}$ and 
 $S = M/({\mathfrak P}\cap {\mathcal O}(X)^{{\mathcal I}_{H}({\mathfrak P})})$.  
 The fact  (\ref{eqn6-2}) implies
 $M = {\mathcal O}(X)$,  $S = {\mathcal O}(X)/ {\mathfrak P}$ and
 $$ {\rm e}({\mathfrak P}, {\mathfrak P}\cap {\mathcal O}(X)^{{\mathcal I}_H({\mathfrak P})})
= 1.$$
 Then from Proposition \ref{keylemma} and (\ref{eqn6-3})
we infer that 
 \begin{equation*}\begin{aligned} & {\rm e}({\mathfrak P}, {\mathfrak P}\cap {\mathcal O}(X)^{{\mathcal I}_H({\mathfrak P})})
\cdot \sharp_p (\langle {\mathfrak X}(G^0)^M\rangle /\langle{\mathfrak X}(G^0)^{S}\rangle)\\
& \qquad= \sharp_p (\langle {\mathfrak X}(G^0)^M\rangle /\langle{\mathfrak X}(G^0)^{S}\rangle)\\
& \qquad \leqq  {\rm e} ({\mathfrak P}, {\mathfrak P}\cap {\mathcal O}(X)^{G^{0}})\\
& \qquad < {\rm e} ({\mathfrak P}, {\mathfrak P}\cap {\mathcal O}(X)^{G^{0}})\cdot p = {\rm e}({\mathfrak P}, {\mathfrak P}\cap {\mathcal O}(X)^{H}). 
\end{aligned} \end{equation*}
This conflicts with the equality of  {\it (ii)} for $X$, $H$ and 
${\mathfrak P}$.

 {\it (i)} $\Rightarrow$ {\it (ii)}:  We may suppose that $G = H  = Z_H(G^0)$ and 
 preserve notations in the circumstances as  in {\it (ii)}, i.e., let
 $(X, G)$ be an arbitrary  {\it effective} regular action 
 of $G$ on any affine Krull $K$-scheme $X = {\rm Spec} (R)$ such that  ${\mathcal Q}(R)^{G^{0}} ={ \mathcal Q}(R^{G^{0}})$ and 
 ${\mathfrak P}\in
{\rm Ht}_1(R, R^G)$.  For any finite subgroup $L$
 of $G = Z_{G}(G^{0})$, by Lemma \ref{qf} 
 we see  $${\mathcal Q}(R)^{L \cdot G^{0}} 
 ={ \mathcal Q}(R^{L\cdot G^{0}}).$$  
By  (\ref{eq1}) and Proposition \ref{mainofpreviouspaper},  we have the equality
  $${\rm e}_{p'}({\mathfrak P}, {\mathfrak P}\cap R^G) =
  \sharp_{p'} ({\mathcal I}_{G}({\mathfrak P})
  \vert_{R}) = {\rm e}_{p'}({\mathfrak P}, {\mathfrak P}\cap R^{{\mathcal I}_G({\mathfrak P})})$$
  of $p'$-parts of ramification indices.    
 Consequently, in order to prove {\it (ii)}, we suppose that $p >0$ and  it suffices  only to show the following equality 
\begin{equation}\label{eqn6-20}\begin{split} {\rm e}_{p}({\mathfrak P}, {\mathfrak P}\cap R^{{\mathcal I}_G({\mathfrak P})})\cdot &
\displaystyle{\sharp_{p} \left(\langle {\mathfrak X}(G^0)^{R^{{\mathcal I}_{G}({\mathfrak P})}}\rangle /\langle{\mathfrak X}(G^0)^{R^{{\mathcal I}_{G}({\mathfrak P})}/{\mathfrak P}\cap R^{{\mathcal I}_{G}({\mathfrak P})}}\rangle\right)}
\\ & ={\rm e}_{p}({\mathfrak P}, {\mathfrak P}\cap R^G).\end{split}  \end{equation}

Put  ${\mathfrak p} = {\mathfrak P}\cap R^{G^0}$ and  $U =\langle {\rm Unip}({\mathcal I}_G({\mathfrak p}))\rangle$.  By Proposition \ref{lemmaunip},  the group $U$
is  finite.  Since ${\mathcal I}_{G}({\mathfrak p})  \supseteq {\mathcal I}_{G}({\mathfrak P})$
and ${\mathcal I}_{G}({\mathfrak p}) \unrhd U$,  by the definition of $U$ we see
$${\mathcal I}_{G}({\mathfrak P}) \unrhd  U \cap {\mathcal I}_{G}({\mathfrak P}) 
= {\mathcal I}_{U}({\mathfrak P})$$
and  the factor group $$ {\mathcal I}_G({\mathfrak P})/{\mathcal I}_U({\mathfrak P})$$
is of order which is not divisible by $p$. Thus  
$$\left(R^{{\mathcal I}_U({\mathfrak P})}/{\mathfrak P}\cap R^{{\mathcal I}_U({\mathfrak P})}\right)
^{{\mathcal I}_G({\mathfrak P})/{\mathcal I}_U({\mathfrak P})}
= R^{{\mathcal I}_G({\mathfrak P})}/{\mathfrak P}\cap R^{{\mathcal I}_G({\mathfrak P})}$$ and 
by Proposition \ref{finitequotientcharactergroup}, we  see 
\begin{equation}\label{eqn6-16}\sharp_{p}\left( \langle {\mathfrak X}(G^0)^{R^{{\mathcal I}_U({\mathfrak P})}/{\mathfrak P}\cap R^{{\mathcal I}_U({\mathfrak P})}}\rangle
/\langle {\mathfrak X}(G^0)^{R^{{\mathcal I}_G({\mathfrak P})}/{\mathfrak P}\cap R^{{\mathcal I}_G({\mathfrak P})}}
\rangle\right) =1.\end{equation}
For the same reason, 
\begin{equation}\label{eqn6-16-1}\sharp_{p}\left( \langle {\mathfrak X}(G^0)^{R^{{\mathcal I}_U({\mathfrak P})}}\rangle
/\langle {\mathfrak X}(G^0)^{R^{{\mathcal I}_G({\mathfrak P})}}
\rangle\right) =1.\end{equation}
By Proposition \ref{key} we  must have
$$\sharp_{p}\left( \langle {\mathfrak X}(G^0)^{R^{{\mathcal I}_U({\mathfrak p})}/{\mathfrak P}\cap R^{{\mathcal I}_U({\mathfrak P})}}\rangle
/\langle {\mathfrak X}(G^0)^{R^{U}/{\mathfrak P}\cap R^{U}}\rangle\right) =1$$
and
\begin{equation}\label{eqn6-17}\begin{split}&  \displaystyle{\sharp_{p} \left(\langle{\mathfrak X}(G^0)^{R^{U}}\rangle  / \langle{\mathfrak X}  (G^0)^{R^{U}/{\mathfrak P}\cap R^{U}}\rangle\right)} \\
&= \sharp_{p}\left(\langle{\mathfrak X}(G^0)^{R^{{\mathcal I}_U({\mathfrak P})}}\rangle / \langle {\mathfrak X}(G^0)^{R^{{\mathcal I}_U({\mathfrak P})}/{\mathfrak P}\cap R^{{\mathcal I}_U({\mathfrak P})}}\rangle
\right)
\end{split}\end{equation}
\begin{equation} = \sharp_{p}\left(\langle {\mathfrak X}(G^0)^{R^{{\mathcal I}_G({\mathfrak P})}}\rangle /\langle {\mathfrak X}(G^0)^{R^{{\mathcal I}_G({\mathfrak P})}/{\mathfrak P}\cap R^{{\mathcal I}_G({\mathfrak P})}}\rangle \right)\label{eqn6-18}\end{equation}
whose last equality follows from (\ref{eqn6-16}) and (\ref{eqn6-16-1}). 
On the other hand, applying  Proposition \ref{keylemma} to the diagonal action
$G^{0}$ on $R^{U}$and on $R^{{\mathcal I}_{U}({\mathfrak P})}$ respectively, we also must have 
\begin{equation}\label{eqn6-19} \begin{aligned} &{\rm e}({\mathfrak P}\cap R^{U}, {\mathfrak P}\cap R^{U\cdot G^0}) =  \displaystyle{\sharp \left(\langle {\mathfrak X}(G^0)^{R^{U}}\rangle /\langle{\mathfrak X}(G^0)^{R^{U}/{\mathfrak P}\cap R^{U}}\rangle\right)}\\
 &{\rm e}({\mathfrak P}\cap R^{{\mathcal I}_{U}({\mathfrak P})}, 
{\mathfrak P}\cap R^{{\mathcal I}_{U}({\mathfrak P})\cdot G^0}) =  \displaystyle{\sharp \left(\langle {\mathfrak X}(G^0)^{R^{{\mathcal I}_{U}({\mathfrak P})}}\rangle /\langle{\mathfrak X}(G^0)^{R^{{\mathcal I}_{U}({\mathfrak P})}/{\mathfrak P}\cap R^{{\mathcal I}_{U}({\mathfrak P})}}\rangle\right)}\end{aligned} \end{equation}
These equalities and (\ref{eqn6-17}) imply
\begin{equation}\label{eqn6-21} {\rm e}_{p}({\mathfrak P}\cap R^{U}, {\mathfrak P}\cap R^{U\cdot G^0})
=  {\rm e}_{p}({\mathfrak P}\cap R^{{\mathcal I}_{U}({\mathfrak P})}, 
{\mathfrak P}\cap R^{{\mathcal I}_{U}({\mathfrak P})\cdot G^0}).\end{equation}
By {\it (iii)} of  Proposition \ref{pr2-7}, 
\begin{equation}\label{eqn6-22} {\rm e}_{p}({\mathfrak P}\cap R^{U\cdot G^0}, {\mathfrak P}\cap R^{G})~ (= {\rm e}_{p}({\mathfrak p}\cap R^{U\cdot G^0}, {\mathfrak p}\cap R^{G}) )=1.\end{equation}
Clearly ${\rm e}({\mathfrak P}, {\mathfrak P}\cap R^{{\mathcal I}_{U}({\mathfrak P})})
={\rm e}({\mathfrak P}, {\mathfrak P}\cap R^{U})$. 
Consider the action ${\mathcal I}_{G}({\mathfrak P})/{\mathcal I}_{U}({\mathfrak P})$
on $R^{{{\mathcal I}_U({\mathfrak P})}}$. Because $${\mathcal I}_{{\mathcal I}_{G}({\mathfrak P})/{\mathcal I}_{U}({\mathfrak P})} ({\mathfrak P}\cap R^{{\mathcal I}_{U} ({\mathfrak P})})
= {\mathcal I}_{G}({\mathfrak P})/{\mathcal I}_{U}({\mathfrak P}), $$ by  ramification theory  
we see that the ramification index
 $${\rm e}({\mathfrak P}\cap R^{{\mathcal I}_U({\mathfrak P})}, {\mathfrak P}\cap R^{{\mathcal I}_G({\mathfrak P})})$$ is a divisor of $\displaystyle\sharp\left( {\mathcal I}_{G}({\mathfrak P})/{\mathcal I}_{U}({\mathfrak P})
 \right),$ which is not divisible by $p$.  Consequently by
 (\ref{eqn6-17}),  (\ref{eqn6-18}), (\ref{eqn6-19}),  (\ref{eqn6-21}) and (\ref{eqn6-22}) we must have 
\begin{align*} {\rm e}_{p}({\mathfrak P}, {\mathfrak P}\cap R^{G}) &=
{\rm e}_{p}({\mathfrak P}, {\mathfrak P}\cap R^{U})\cdot {\rm e}_{p}({\mathfrak P}\cap R^{U}, {\mathfrak P}\cap R^{U\cdot G^0})\\
&= {\rm e}_{p}({\mathfrak P}, {\mathfrak P}\cap R^{{\mathcal I}_{U}({\mathfrak P})})
\cdot  {\rm e}_{p}({\mathfrak P}\cap R^{{\mathcal I}_{U}({\mathfrak P})}, 
{\mathfrak P}\cap R^{{\mathcal I}_{U}({\mathfrak P})\cdot G^0})\\
&=  {\rm e}_{p}({\mathfrak P}, {\mathfrak P}\cap R^{{\mathcal I}_{U}({\mathfrak P})})\cdot \displaystyle{\sharp_{p} \left(\langle {\mathfrak X}(G^0)^{R^{{\mathcal I}_{U}({\mathfrak P})}}\rangle /\langle{\mathfrak X}(G^0)^{R^{{\mathcal I}_{U}({\mathfrak P})}/{\mathfrak P}\cap R^{{\mathcal I}_{U}({\mathfrak P})}}\rangle\right)}\\
&= {\rm e}_{p}({\mathfrak P}, {\mathfrak P}\cap R^{{\mathcal I}_G({\mathfrak P})})\cdot
\displaystyle{\sharp_{p} \left(\langle {\mathfrak X}(G^0)^{R^{{\mathcal I}_{G}({\mathfrak P})}}\rangle /\langle{\mathfrak X}(G^0)^{R^{{\mathcal I}_{G}({\mathfrak P})}/{\mathfrak P}\cap R^{{\mathcal I}_{G}({\mathfrak P})}}\rangle\right)} \end{align*}
which shows (\ref{eqn6-20})  required as above. Thus the proof is completed. \qed
  
  \begin{exam} Case `` $G^{0}$ is a simple algebraic group'' : \rm
  Let $k$ be an {\it odd} natural number $> 1$
  such that $k$ is not divisible by $p$ if $p>0$ and $\zeta_{k}$ be 
  a fixed primitive $k$-th root of $1 \in K$.  Define the subgroup 
  
  $$G := \left\langle  SL_{2}(K), \left[ \begin{array}{cc} 
  \zeta_{k} & 0 \\ 0 & \zeta_{k}  \end{array} \right] \right\rangle \subseteq
  GL_{2}(K)$$ 
  which acts naturally on the $K$-vector space $KX_{1}\oplus KX_{2}$. 
  The group $G$ is an affine algebraic group with  $G^{0} = SL_{2}(K)$ and $G = Z_{G}(G^{0}) \cdot G^{0}$. 
  Let $R$ be a polynomial ring $K[Y_{1}, Y_{2}, Y_{3}]$ over $K$
  with three variables.   Identifying
    $$Y_{1} = X_{1}^{2}, 
  Y_{2}  = X_{1}X_{2} ~~\text{and} ~~
  Y_{3} = X_{2}^{2},$$  
the group  $G$ acts on  $R$ induced from  the representation
  on the  homogeneous
  part of $K[X_{1}, X_{2}]$ of degree $2$. Then there is a 
  homogeneous polynomial $f$ of $R$ of degree $2$  such that
  $R^{G^{0}} = K[f]$.  Clearly $f$ is a prime element of $R$ and ${\mathfrak P}: = Rf \in {\rm Ht}_{1}(R, R^{G})$. 
  For any finite subgroup $\Gamma$ of $G$, we have 
  $${\rm e}( {\mathfrak P}, {\mathfrak P}\cap R^{\Gamma})
  = {\rm e}( {\mathfrak P}, {\mathfrak P}\cap R^{{\mathcal I}_{\Gamma}({\mathfrak P})})
   =1,$$
   because $${\mathcal I}_{\Gamma}({\mathfrak P}) \subseteq 
   {\mathfrak R}(R, G) \subseteq Z_{G}(G^{0}) \subseteq Z(GL_{2}(K))~~ (cf. ~\cite{Nak5}).$$
  However
  $${\rm e}( {\mathfrak P}, {\mathfrak P}\cap R^{G})  = k > {\rm e}( {\mathfrak P}, {\mathfrak P}\cap R^{\Gamma}),$$
  which shows that  {\it Problem 1.1 is not affirmative} in this case. \end{exam}
  
 
 \medskip 
 

\end{document}